\pdfoutput=1
\documentclass{article}

\usepackage{iclr2027_conference,times}

\usepackage[utf8]{inputenc}
\usepackage[T1]{fontenc}
\usepackage{microtype}

\usepackage{amsmath, amssymb, amsthm, mathtools}
\usepackage{bm}
\usepackage{graphicx}
\usepackage[colorlinks=true, allcolors=blue]{hyperref}
\usepackage{url}
\usepackage{booktabs}
\usepackage{algorithm}
\usepackage{algpseudocode}
\usepackage{enumitem}
\usepackage{tikz}
\usetikzlibrary{positioning, arrows.meta, fit, calc}
\usepackage{multirow}
\usepackage{subcaption}

\newtheorem{theorem}{Theorem}[section]
\newtheorem{proposition}[theorem]{Proposition}

\newtheorem{corollary}[theorem]{Corollary}
\theoremstyle{definition}

\theoremstyle{remark}

\newcommand{\R}{\mathbb{R}}
\newcommand{\E}{\mathbb{E}}
\newcommand{\PP}{\mathbb{P}}
\newcommand{\argmin}{\operatorname*{arg\,min}}

\newcommand{\vect}[1]{\mathbf{#1}}
\newcommand{\ba}{\vect{a}}
\newcommand{\bs}{\vect{s}}
\newcommand{\bx}{\vect{x}}
\newcommand{\by}{\vect{y}}

\newcommand{\bo}{\vect{o}}
\newcommand{\bu}{\vect{u}}

\newcommand{\blam}{\bm{\lambda}}
\newcommand{\bmu}{\bm{\mu}}
\newcommand{\bth}{\bm{\theta}}
\newcommand{\bpsi}{\bm{\psi}}

\newcommand{\bpi}{\bm{\pi}}

\title{Learning Advantage Function Surrogates \\ via Convex Conjugate Matching}
\author{Alvaro Maggiar \\
  Amazon.com \\
  \texttt{maggiara@amazon.com}
}

\iclrfinalcopy

\begin{document}
\maketitle
\lhead{}

\begin{abstract}
A policy trained to make per-unit decisions, such as how much of a product to order, tells us the best action but not the cost of deviating from it.  That cost, the advantage function, is what downstream systems need when they have to depart from the per-unit optimum: to estimate the impact of ordering more or fewer units, to respond to a one-off vendor discount, or to satisfy a shared capacity constraint. The advantage function is what makes such trade-offs computable, yet pure policy-learning algorithms often produce the policy alone, without a critic that would yield it. We propose a method to train an advantage function surrogate through price perturbations motivated by a Fenchel-Young loss. We describe a two-stage recipe (fit the policy, then the cost model) and illustrate it on a multi-period inventory management problem.
\end{abstract}

\section{Introduction}
\label{sec:intro}

\paragraph{Motivation.}
Consider a policy trained to take optimal decisions at the level of individual products, accounts, or units, such as how many units of a given SKU to buy, learned in a setting where each unit targets its per-unit optimum. Beyond the optimal action, several downstream uses call for an \emph{advantage} function that quantifies the cost of deviating from this optimum: counterfactual impact assessment (the financial effect of buying $X$ more or fewer units), one-off vendor offers (the order quantity trades a time-limited discount against the advantage's curvature around the optimum), and occasional coupling constraints (vendor minimum-order quantities or capacity limits requiring coordinated deviations across products). In each case the policy is unchanged; what is needed is a faithful, typically convex, model of the cost landscape around the optimum.

\paragraph{A generic consumption model.}
For a unified abstraction of the uses above, suppose an external consumer observes state $\bs$ and solves
\begin{equation}
\label{eq:downstream}
\ba^{\mathcal{E}}(\bs) \in \argmin_{\ba } \big\{\, Q(\ba; \bs) + g(\ba; \bs) \,\big\},
\end{equation}
where $g$ is a convex objective specific to the consumer and $Q$ is the action-value function. Setting $g \equiv 0$ recovers the unconstrained optimum, a linear $g(\ba) = -\blam^\top \ba$ models a vendor discount or the dualization of a coupling constraint, and a quadratic or coupled $g$ captures more elaborate business rules. We aim for a surrogate of $Q$, and this surrogate need not approximate $Q$ pointwise, what matters is to produce the right decision under \eqref{eq:downstream}. In particular, when $Q(\cdot;\bs)$ is strictly convex and $g$ is differentiable, $\ba^{\mathcal{E}}(\bs) = \nabla Q^*(-\nabla g(\ba^{\mathcal{E}};\bs);\bs)$, where $Q^*$ is the (state-conditional) convex conjugate of $Q$ in the action argument. A good surrogate $f_{\bmu}$ is thus one such that $\nabla f^*_{\bmu}(\blam;\bs) \approx \nabla Q^*(\blam;\bs)$ for dual values $\blam$ typically encountered in the consumption uses, so that the consumer's actions under the surrogate match those under the true action-value function.

\paragraph{Why the policy alone is not enough.}
Many policy-learning methods return $\bpi_{\bth}$ without $Q$, including the direct policy optimization through differentiable simulators \citep{madeka2022deep,alvo2023neural} that we use upstream. The standard remedy fits a critic by regression on Bellman or Monte-Carlo targets (Section~\ref{sec:related}). A critic fit by pointwise regression is accurate where its data lie, which for data collected under a near-deterministic policy is close to the optimum, while the consumer queries $Q$ at the actions $\nabla Q^*(\blam;\bs)$ to which prices push it. A second option learns the surrogate jointly with the policy, through a differentiable optimization layer \citep{amos2017optnet,agrawal2019differentiable} whose argmin is the action. This is under-determined, since any surrogate whose minimum coincides with the optimal action yields the same reward: for a quadratic $A(\ba;\bs)=\tfrac{1}{2}(\ba-\ba^*(\bs))^\top B (\ba-\ba^*(\bs))$, every positive-definite $B$ is indistinguishable so long as $\ba^*(\bs)$ is optimal.

\paragraph{Approach.}
To resolve the under-determinism we fit a strongly convex advantage surrogate $f_{\bmu}(\ba; \bs)$ around the optimal action under a Fenchel-Young loss, trained against a distribution $p(\blam)$ of dual perturbations that simulate the shadow prices the consumption uses induce. The distribution $p(\blam)$ is an \emph{exogenous} design choice: the policy remains a function of $\bs$ alone, optimal in the unconstrained problem, and the role of $p(\blam)$ is to span the range of dual values the surrogate will be queried at. We develop a two-stage procedure (policy first, surrogate second). The method applies to any strongly convex family whose conjugate gradient can be evaluated and differentiated, and the family can encode problem knowledge, for example a quadratic when the consumer solves QPs, or discrete-convex structure such as $M^\natural$-convexity \citep{murota2003discrete}, regularized to be strongly convex.

\paragraph{Contributions.}
\begin{itemize}[leftmargin=*, itemsep=1pt, topsep=2pt]
\item We cast the consumption of an advantage function as a price-perturbed argmin \eqref{eq:downstream}, and show that fitting a surrogate through its unperturbed argmin leaves its curvature unidentified (Section~\ref{subsec:underid}).
\item We propose a conjugate-matching objective derived from the Fenchel-Young inequality. It requires only evaluations of $Q^{\bth}$ at the surrogate's own decisions, never the consumer's optimal decisions, and, for convex $Q^{\bth}$, we bound the error of the decisions it induces and characterize what it identifies (Proposition~\ref{prop:decerr}, Corollary~\ref{prop:ident}) and its stationary points (Theorems~\ref{thm:foc}--\ref{thm:pge}).
\item We give an anchored parameterization and a two-stage protocol that fit the surrogate to a frozen, possibly deployed, policy, and illustrate them on a lost-sales inventory problem with a shared capacity constraint represented by its shadow price.
\end{itemize}
Section~\ref{sec:related} positions the method against related work, Section~\ref{sec:prelim} sets up the exo-MDP and convex-analytic background, Section~\ref{sec:method} develops the surrogate and its loss, Section~\ref{sec:algorithm} the training procedure and its requirements, and Section~\ref{sec:experiments} the experiments.

\section{Related Work}
\label{sec:related}

\paragraph{Critics with tractable action structure.}
Actor-critic methods learn $Q$ or $A$ by regression on Bellman or Monte-Carlo targets \citep{sutton1999policy,silver2014deterministic,haarnoja2018soft}. Normalized Advantage Functions \citep{gu2016continuous} restrict the advantage to a quadratic $-\tfrac12(\ba-\mu(\bs))^\top P(\bs)(\ba-\mu(\bs))$ so that the greedy action is closed-form, and input convex neural networks \citep{amos2017input} make $Q$ convex in the action so that the greedy action solves a convex program. In both, the critic is fit pointwise by Bellman regression and the policy is the critic's own greedy action. Our surrogate class contains the NAF quadratic, and an ICNN once a strongly convex term is added, but it is fit after the policy, to a frozen $\bpi_{\bth}$, and only through the decisions it induces under price perturbations, so its accuracy is concentrated on the prices at which consumers query it.

\paragraph{Fenchel-Young losses and perturbed optimizers.}
\citet{blondel2020learning} define the Fenchel-Young loss $\Omega(\boldsymbol{v})+\Omega^*(\boldsymbol{\eta})-\langle\boldsymbol{\eta},\boldsymbol{v}\rangle$ for a fixed regularizer $\Omega$, a predicted dual argument $\boldsymbol{\eta}$ and an observed target $\boldsymbol{v}$. \citet{berthet2020learning} stochastically perturb the input of a discrete optimizer with a linear objective, which makes an otherwise locally constant argmax differentiable with smoothness set by the noise amplitude, and connect the resulting perturbed optimizers to Fenchel-Young losses. Our loss \eqref{eq:fenchel_young_surrogate} uses the same inequality with the roles rearranged: the function is the unknown $Q^{\bth}$, which need not be convex and is available only through evaluation; the dual argument $\blam$ is drawn from a distribution we design; and the predicted object is the primal action $\nabla f^*_{\bmu}(\blam;\bs)$. Since the surrogate is strongly convex and already differentiable, the perturbation plays no smoothing role here; its purpose is to identify the surrogate's curvature.

\paragraph{Inverse optimization and decision-focused surrogates.}
Inverse optimization estimates the parameters of an optimization model from observed, possibly noisy, solutions \citep{aswani2018inverse}, and \citet{li2025fenchel} add a strongly convex regularizer to a linear forward problem and show that the suboptimality loss of an observed decision then coincides with a Fenchel-Young loss. \citet{gupta2024data} learn a simpler, for example convex, surrogate optimization model by minimizing the squared distance between its optimal solutions and those of the original model, on a dataset obtained by solving the original model offline for sampled inputs. Both learn from decision data: observed, possibly noisy or suboptimal, decisions in inverse optimization, and optimal solutions of the original model in decision-focused surrogate modeling. In our setting the decisions $\argmin_{\ba}\{Q^{\bth}(\bs,\ba)-\blam^\top\ba\}$ that the consumer would take under a price are not available, and computing each would require a separate optimization over rollouts; the loss \eqref{eq:training-objective} only evaluates $Q^{\bth}$ at the surrogate's own decision, which a differentiable rollout provides, and measures the error in cost rather than in distance between decisions. In this respect the objective is closer to decision-focused learning \citep{donti2017task,wilder2019melding,elmachtoub2022smart}, which trains a model of the parameters of an optimization problem on the cost of the decisions it induces, or on a convex surrogate of that cost. Here the learned object is the objective itself, and the randomness is over the price the consumer imposes.

\paragraph{Conditioning the policy on the constraint.}
For a known constraint family, one can instead train the policy, or an actor-critic, with the constraint parameter as an input. When the constraint couples many products, the relevant input after dualization is the shared shadow price $\blam$, and policies taking it as an input have been trained \citep{eisenach2024neural} and coordinated through a broadcast price \citep{eisenach2024neural,wang2026ready}. Such a price-indexed policy $\pi(\bs,\blam)$ is also the natural object for coordination through a dual interface \citep{maggiar2026consensus}. The surrogate induces one, $\nabla f^*_{\bmu}(\blam;\bs)$, with three differences. First, the price must enter the training of a price-indexed policy, so a policy already trained without it cannot be used, whereas the surrogate is fit to a frozen policy. Second, a price-indexed policy is consistent with a convex cost only if it is the gradient of a convex function of $\blam$, which requires input-convex or gradient-of-convex architectures; the surrogate has this property by construction. Third, with only a price-to-action map a consumer with a nonlinear objective $g$ must solve the fixed point $\ba = \pi(\bs, -\nabla g(\ba))$, whereas with the cost $f_{\bmu}$ it minimizes $f_{\bmu} + g$ directly; the dual function of $f_{\bmu}$ is also differentiable, so a shared price can be set by dual gradient ascent.

\section{Preliminaries}
\label{sec:prelim}

\subsection{Exogenous MDPs}
\label{subsec:exomdp}

We work with exogenous MDPs (exo-MDPs) \citep{sinclair2023hindsight}, a subclass of the input-driven MDPs of \citet{mao2018variance}, in which uncertainty is driven by an exogenous process independent of the agent's actions. At each epoch $t \in \{0, 1, \ldots, T-1\}$, the state decomposes as $\bs_t = (\bx_t, \by_t)$, where $\bx_t = (\bo_{t-H}, \ldots, \bo_{t-1})$ is a fixed-length history of exogenous observations \citep{madeka2022deep}, and $\by_t$ the endogenous part of the state. The agent's action $\ba_t$ is selected by a deterministic policy $\bpi_{\bth}$ parameterized by $\bth$, $\ba_t = \bpi_{\bth}(\bs_t)$. The endogenous state then evolves under a deterministic, differentiable transition $\by_{t+1} = \mathcal{T}(\bs_t, \bo_t, \ba_t)$, and a differentiable per-period cost $c_t(\bs_t, \ba_t)$ is incurred. The exogenous sequence $\{\bo_t\}_t$ is drawn from a distribution $\PP$ independent of the agent.

The objective is to minimize the finite-horizon discounted expected cost, or the long-run average cost in the stationary case,
\begin{align*}
J_T(\bth; \by_0) = \E_{\PP}\!\left[\sum_{t=0}^{T-1} \gamma^t\, c_t(\bs_t, \ba_t)\,\Big|\,\by_0\right],
\qquad
J(\bth) = \lim_{T\to\infty} \E_{\PP}\!\left[\tfrac{1}{T}\sum_{t=0}^{T-1}\, c(\bs_t, \ba_t)\right].
\end{align*}
The exogenous structure lets objectives and gradients be evaluated by rolling out against stored historical traces, which underlies the DirectBackprop/HDPO training algorithms \citep{madeka2022deep, alvo2023neural} that we rely on for upstream training.

\subsection{Policy and Value Functions}
\label{subsec:upstream}

We assume that a trained instance of the upstream policy from Section~\ref{subsec:exomdp} is available, namely $\bpi_{\bth}$ with $\bth$ fixed to its trained value. The value, action-value, and advantage functions of $\bpi_{\bth}$ (in cost form) are
\[
V_t^{\bth}(\bs) \;=\; \E\!\left[\sum_{\tau = t}^{T-1} \gamma^{\tau - t}\, c_\tau(\bs_\tau, \ba_\tau) \,\Big|\, \bs_t = \bs,\ \ba_\tau = \bpi_{\bth}(\bs_\tau)\right],
\]
\[
Q_t^{\bth}(\bs, \ba) \;=\; c(\bs, \ba) + \gamma\, \E\!\left[V_{t+1}^{\bth}(\bs_{t+1}) \,\big|\, \bs_t = \bs,\ \ba_t = \ba\right],
\qquad
A_t^{\bth}(\bs, \ba) \;=\; Q_t^{\bth}(\bs, \ba) - V_t^{\bth}(\bs).
\]
In the infinite-horizon stationary setting (time-homogeneous dynamics and cost), these are independent of $t$ and we drop the subscript, writing $V^{\bth}$, $Q^{\bth}$, $A^{\bth}$. We adopt this stationary convention for the remainder of the paper for ease of exposition.

By definition $A^{\bth}(\bs, \bpi_{\bth}(\bs)) = 0$; if $\bpi_{\bth}$ is moreover greedy with respect to $Q^{\bth}$, as an optimal policy is, then $A^{\bth}(\bs, \ba) \geq 0$ for all $\ba$.

When the exo-MDP decomposition is not relevant, we write $\bs$ for the full state.

\subsection{Convex-Analytic Background}
\label{subsec:convex}

We use several standard facts related to the convex conjugate (\citealp[Sections~12 and~26]{rockafellar1997convex}; \citealp[Prop.~2]{blondel2020learning}). For the convex conjugate $f^*(\blam) := \sup_{\ba}\{\langle \ba,\blam\rangle - f(\ba)\}$, if $f$ is closed and strictly convex, $f^*$ is differentiable on $\operatorname{int}\operatorname{dom} f^*$, where the perturbed minimizer $\argmin_{\ba}\{f(\ba) - \blam^\top \ba\}$ exists and equals $\nabla f^*(\blam)$; for differentiable $f$ this is the Legendre correspondence $\ba = \nabla f^*(\blam) \Leftrightarrow \blam = \nabla f(\ba)$. The Fenchel-Young gap $L^{\mathrm{FY}}(\ba, \blam; f) := f(\ba) + f^*(\blam) - \langle \ba, \blam\rangle \geq 0$ vanishes iff $\blam \in \partial f(\ba)$ ($\blam = \nabla f(\ba)$ for differentiable $f$) and can serve as a fitting loss \citep{blondel2020learning}. For $C^2$ parameterized $f_{\bmu}$ with $\nabla^2_{\ba} f_{\bmu} \succ 0$, implicit differentiation of $\nabla_{\ba} f_{\bmu}(\ba) = \blam$ gives the argmin Jacobian $\tfrac{\partial \ba}{\partial \bmu} = -[\nabla_{\ba}^2 f_{\bmu}(\ba)]^{-1}\nabla_{\bmu}\nabla_{\ba} f_{\bmu}(\ba)$, the OptNet KKT Jacobian \citep{amos2017optnet} for quadratic $f_{\bmu}$ and the layer of \citet{agrawal2019differentiable} for DPP-representable families.

\section{The Advantage Surrogate via Convex Conjugate Matching}
\label{sec:method}

\subsection{Surrogate Class}
\label{subsec:surrogate-class}

\paragraph{General parametric family.} Let $\mathcal{F} = \{f_{\bmu}(\cdot; \cdot) : \bmu \in \mathcal{M}\}$ be a family of functions $f_{\bmu} : \R^d \times \mathcal{S} \to \R$, parameterized by $\bmu \in \mathcal{M}$. We require $f_{\bmu}(\cdot; \bs)$ to be strongly convex in $\ba$ with continuous gradient $\nabla_{\ba} f_{\bmu}$, and with strong-convexity modulus bounded below uniformly on bounded subsets of $\mathcal{S}$. Under these conditions the conjugate $f_{\bmu}^*(\blam; \bs) = \sup_{\ba}\{\langle \ba, \blam\rangle - f_{\bmu}(\ba; \bs)\}$ is finite everywhere and differentiable in $\blam$ on $\R^d$.

\paragraph{Anchored parameterization.} A useful specialization factorizes $f$ into a \emph{decision} component and a \emph{shape} component:
\begin{equation}
\label{eq:anchored}
f_{\bmu, \bth}(\ba; \bs) = \psi_{\bmu}\!\left(\ba - \ba^*(\bs; \bth);\; \bs\right),
\end{equation}
where $\ba^*(\cdot; \bth) : \mathcal{S} \to \R^d$ is the \emph{anchor network} representing the optimal policy, and $\psi_{\bmu} : \R^d \times \mathcal{S} \to \R$ is a \emph{shape function} of the deviation $\bu = \ba - \ba^*(\bs;\bth)$, strongly convex with continuous gradient in its first argument (we use $C^2$ when convenient), satisfying
\[
\psi_{\bmu}(\mathbf{0}; \bs) = 0, \qquad \nabla_{\bu} \psi_{\bmu}(\mathbf{0}; \bs) = \mathbf{0} \qquad \text{for every } (\bmu, \bs).
\]
These conditions enforce $f_{\bmu, \bth}(\ba^*; \bs) = 0$ and $\nabla_{\ba} f_{\bmu, \bth}(\ba^*; \bs) = \mathbf{0}$ structurally, and the conjugate separates cleanly as $f_{\bmu, \bth}^*(\blam; \bs) = \langle \ba^*(\bs; \bth), \blam\rangle + \psi_{\bmu}^*(\blam; \bs)$, so the price-perturbed argmin is $\nabla_{\blam} f_{\bmu, \bth}^*(\blam; \bs) = \ba^*(\bs; \bth) + \nabla_{\blam}\psi_{\bmu}^*(\blam; \bs)$. The factorization supports \emph{two-stage} training (anchor first, shape second; Section~\ref{sec:algorithm}).

\paragraph{Example (quadratic).} Taking $\psi_{\bmu}(\bu; \bs) = \tfrac{1}{2}\bu^\top B(\bs; \bmu)\bu$ with $B \succ 0$ recovers the anchored quadratic with closed-form conjugate $\psi_{\bmu}^*(\blam; \bs) = \tfrac{1}{2}\blam^\top B(\bs; \bmu)^{-1}\blam$. The same algebraic form, with $(\mu(\bs), P(\bs))$ in place of $(\ba^*, B)$, is the Normalized Advantage Function $Q$-network of \citet{gu2016continuous}, which differs in its loss and protocol (Section~\ref{sec:related}).

\subsection{Under-Identification of the Simple Loss}
\label{subsec:underid}

Directly minimizing the upstream cumulative cost with the differentiable-argmin action $\ba^{\star}_{\bmu}(\bs) := \argmin_{\ba} f_{\bmu}(\ba; \bs)$ fails to identify the surrogate beyond its argmin. The naive loss $\mathcal{L}_{\mathrm{naive}}(\bmu) := \E[\sum_t \gamma^t c(\bs_t, \ba^{\star}_{\bmu}(\bs_t))]$, with the trajectory rolled out under $\ba_t = \ba^{\star}_{\bmu}(\bs_t)$, depends on $\bmu$ only through the argmin map: if $\ba^{\star}_{\bmu} = \ba^{\star}_{\bmu'}$ on every state the rollout can reach then $\mathcal{L}_{\mathrm{naive}}(\bmu) = \mathcal{L}_{\mathrm{naive}}(\bmu')$. For the anchored parameterization \eqref{eq:anchored} the argmin is $\ba^*(\bs; \bth)$, independent of the curvature $\bmu$, so all curvatures minimize $\mathcal{L}_{\mathrm{naive}}$ equally.

\subsection{Convex Conjugate Matching}

To break this degeneracy, we learn the action-value function through its convex conjugate. Let $Q^{\bth*}(\bs,\blam) = \sup_{\ba}\{\blam^\top \ba - Q^{\bth}(\bs,\ba)\}$ be the state-conditional conjugate in the action argument, and let the surrogate-induced action be the perturbed argmin $\ba(\blam, \bs; \bmu) := \argmin_{\ba}\{f_{\bmu}(\ba; \bs) - \blam^\top \ba\} = \nabla f_{\bmu}^*(\blam; \bs)$, which collapses to the anchor $\ba^{\star}_{\bmu}(\bs)$ at $\blam = \mathbf 0$ and traces the dual conjugate of $f_{\bmu}$ for $\blam \neq \mathbf 0$. The Fenchel-Young inequality applied to this action,
\begin{align}
Q^{\bth}(\bs,\ba(\blam, \bs; \bmu)) + Q^{\bth*}(\bs,\blam) - \blam^\top \ba(\blam, \bs; \bmu) \;\geq\; 0,\label{eq:fenchel_young_surrogate}
\end{align}
holds with equality iff the surrogate's action minimizes $Q^{\bth}(\bs,\cdot)-\blam^\top\ba$, that is, for differentiable convex $Q^{\bth}$, iff $\nabla_{\ba} Q^{\bth}(\bs,\ba(\blam,\bs;\bmu)) = \blam$; when the minimizer is unique, iff $\ba(\blam, \bs; \bmu)=\ba(\blam, \bs):=\argmin_{\ba}\{Q^{\bth}(\bs,\ba)-\blam^\top\ba\}$. A surrogate yielding decisions close to the true action-value function thus minimizes the gap in \eqref{eq:fenchel_young_surrogate}. (If $Q^{\bth}$ is not convex in $\ba$, it has the same conjugate as its closed convex envelope $Q^{\bth**}$, so the loss can at best match $Q^{\bth**}$; the results below assume convexity.)

\subsection{Fenchel-Young Loss (Two-Stage)}
\label{subsec:fy-loss}

Assume the upstream policy $\bpi_{\bth}$ and value function $V^{\bth}$ are given (Section~\ref{subsec:upstream}). Let $\rho^{\bth}$ be the distribution of the state and next exogenous input $(\bs_t,\bo_t)$ along $\bpi_{\bth}$-rollouts and $p(\blam)$ a chosen dual-perturbation distribution. Averaging the FY gap \eqref{eq:fenchel_young_surrogate} over $\rho^{\bth}\otimes p$ and dropping the $\bmu$-independent constant $\E[Q^{\bth*}(\bs,\blam)]$, assumed finite, the surrogate is fit by minimizing
\begin{align}
\mathcal{J}(\bmu)&:=\E_{(\bs_t,\bo_t) \sim \rho^{\bth},\, \blam \sim p(\blam)}\!\left[\, Q^{\bth}(\bs_t, \ba(\blam, \bs_t; \bmu)) \,-\, \blam^\top \ba(\blam, \bs_t; \bmu) \,\right]\label{eq:training-objective}\\
&\phantom{:}=\E_{(\bs_t,\bo_t) \sim \rho^{\bth},\, \blam \sim p(\blam)}\!\left[\, c_t(\bs_t, \ba(\blam, \bs_t; \bmu)) \,-\, \blam^\top \ba(\blam, \bs_t; \bmu) \,+\, \gamma\, V_{t+1}^{\bth}(\bs_{t+1})\,\right],\nonumber
\end{align}
the two forms being equal by the definition of $Q^{\bth}$ and the tower property. The target is exact with a rollout continuation of $V^{\bth}$ and approximate with a learned value network. The states are drawn from the $\bpi_{\bth}$-rollout, keeping the second-stage distribution aligned with deployment.

\paragraph{Decision error and identification.} By \eqref{eq:fenchel_young_surrogate}, the excess loss $\Delta(\bmu) := \mathcal{J}(\bmu) + \E[Q^{\bth*}(\bs,\blam)] \ge 0$ is the average Fenchel-Young gap. The next result relates it to the error in the decisions the surrogate induces (proofs in Appendix~\ref{app:ident}).

\begin{proposition}[Decision error]
\label{prop:decerr}
Assume that for $\rho^{\bth}$-almost every $\bs$, $Q^{\bth}(\bs,\cdot)$ is strictly convex and differentiable on an open set containing a convex set $\mathcal{C}_{\bs}$, the true decision $\ba(\blam,\bs)$ exists for every $\blam \in \operatorname{supp} p$, $\mathcal{C}_{\bs}$ contains $\ba(\blam,\bs)$ and $\ba(\blam,\bs;\bmu)$ for all such $\blam$, and, for constants $0 < m \le L$ independent of $\bs$, $Q^{\bth}(\bs,\cdot)$ is $m$-strongly convex with $L$-Lipschitz gradient on $\mathcal{C}_{\bs}$. Then
\begin{equation*}
\tfrac{m}{2}\,\E\big\|\ba(\blam,\bs;\bmu) - \ba(\blam,\bs)\big\|^2 \;\le\; \Delta(\bmu) \;\le\; \tfrac{L}{2}\,\E\big\|\ba(\blam,\bs;\bmu) - \ba(\blam,\bs)\big\|^2,
\end{equation*}
with expectations under $\rho^{\bth}\otimes p$.
\end{proposition}

Minimizing $\mathcal{J}$ thus minimizes, up to the factor $L/m$, the mean squared error of the surrogate's decisions over the prices in the support of $p$, whether or not the family contains $Q^{\bth}$, and the constants need only hold on the actions these prices reach, for the true and the surrogate's decisions.

\begin{corollary}[Identification]
\label{prop:ident}
Under the assumptions of Proposition~\ref{prop:decerr}, suppose $\Delta(\bmu) = 0$, and for $\rho^{\bth}$-almost every $\bs$ let $\Lambda \subseteq \operatorname{int}\operatorname{dom} Q^{\bth*}(\bs,\cdot)$ be an open convex set on which $p$ has a positive density. Then $f_{\bmu}(\cdot;\bs) - Q^{\bth}(\bs,\cdot)$ is constant on $\mathcal{A}_\Lambda := \nabla Q^{\bth*}(\Lambda;\bs)$. If moreover $\mathbf 0 \in \Lambda$ and $f_{\bmu}$ is anchored at $\ba^*(\bs) = \bpi_{\bth}(\bs)$, then $f_{\bmu}(\cdot;\bs) = A^{\bth}(\bs,\cdot)$ on $\mathcal{A}_\Lambda$.
\end{corollary}

The second statement requires a greedy policy, since $\Delta(\bmu)=0$ and $\mathbf 0 \in \Lambda$ force $\bpi_{\bth}(\bs) \in \argmin Q^{\bth}(\bs,\cdot)$. A point mass $p = \delta_{\mathbf 0}$ identifies only the anchor, which is the under-identification of Section~\ref{subsec:underid}, and the loss depends on $f_{\bmu}$ only through its decisions at prices in the support of $p$, so $p$ must span the prices the consumer will impose.

\paragraph{Example.} Take $d = 1$, a quadratic action-value $Q(a) = \tfrac{k}{2}(a - a^*)^2 + V^{\bth}(\bs)$ with $k > 0$, and the anchored quadratic surrogate $f_{\bmu}(a) = \tfrac{B}{2}(a - a^*)^2$ with $B > 0$. The true and surrogate decisions under the price $\lambda$ are $a_\lambda = a^* + \lambda/k$ and $a_{\bmu} = a^* + \lambda/B$. The curvature of $Q$ is constant, so $m = L = k$ and the bound of Proposition~\ref{prop:decerr} holds with equality:
\[
\Delta(\bmu) = \tfrac{k}{2}\,\E[\lambda^2]\left(\tfrac{1}{B} - \tfrac{1}{k}\right)^2 .
\]
If $\E[\lambda^2] > 0$, the loss vanishes only at $B = k$, where $f_{\bmu} = A^{\bth}$. If $p = \delta_0$, it vanishes for every $B$, which is the under-identification of Section~\ref{subsec:underid}. The second moment $\E[\lambda^2]$ sets the curvature of $\Delta$ in $1/B$, and hence how sharply the loss determines the surrogate's curvature.

\subsection{First-Order Condition and Policy-Gradient Equivalence}
\label{subsec:foc}

Assume $f_{\bmu}$ is $C^2$ with $\nabla^2_{\ba} f_{\bmu} \succ 0$, $Q^{\bth}$ is differentiable in $\ba$, and differentiation and expectation can be exchanged. Differentiating \eqref{eq:training-objective} through the differentiable exo-MDP transition gives $\nabla_{\bmu} \mathcal{J}(\bmu) = \E[\tfrac{\partial \ba}{\partial \bmu}^{\!\top}(\nabla_{\ba} Q^{\bth}(\bs_t,\ba) - \blam)]$. Using $\nabla_{\ba} Q^{\bth} = \nabla_{\ba} A^{\bth}$ ($V^{\bth}$ is action-independent) and the dual consistency $\blam = \nabla_{\ba} f_{\bmu}(\ba; \bs_t)$ at the argmin, the bracket becomes the advantage residual $\nabla_{\ba} A^{\bth} - \nabla_{\ba} f_{\bmu}$.

\begin{theorem}[First-order condition]
\label{thm:foc}
At every stationary point of \eqref{eq:training-objective}, with $\ba = \ba(\blam, \bs; \bmu)$,
\begin{equation}
\label{eq:moment-condition}
\E_{\bs \sim \rho^{\bth},\, \blam \sim p(\blam)}\!\left[\, \frac{\partial \ba}{\partial \bmu}^{\!\top}\,\big(\nabla_{\ba} A^{\bth}(\bs, \ba) - \nabla_{\ba} f_{\bmu}(\ba; \bs)\big)\,\right] = \mathbf{0},
\end{equation}
equivalently with $Q^{\bth}$ in place of $A^{\bth}$.
\end{theorem}

This FOC has a policy-gradient reading. Rearranging \eqref{eq:moment-condition} at $\ba = \ba(\blam, \bs; \bmu)$,
\begin{align*}
\E_{\rho^{\bth} \otimes p}\!\left[\, \tfrac{\partial \ba}{\partial \bmu}^{\!\top}\nabla_{\ba} Q^{\bth}(\bs, \ba) \right] \;=\; \E_{\rho^{\bth} \otimes p}\!\left[\, \tfrac{\partial \ba}{\partial \bmu}^{\!\top}\nabla_{\ba} f_{\bmu}(\ba; \bs) \right],
\end{align*}
each side has the form of an off-policy deterministic policy gradient \citep{silver2014deterministic} through the dual policy, with the upstream $Q^{\bth}$ in place of the dual policy's own action-value on the left and the surrogate $f_{\bmu}$ on the right. The fitting FOC imposes that the two coincide under $\rho^{\bth}\otimes p$.

\begin{theorem}[Policy-gradient equivalence]
\label{thm:pge}
Let $\tilde{\bpi}_{\bmu}(\bs, \blam) := \ba(\blam, \bs; \bmu) = \nabla f_{\bmu}^*(\blam; \bs)$ denote the dual policy induced by $f_{\bmu}$. For any joint distribution $\nu$ on $(\bs, \blam)$ and any cost function $h(\bs, \ba)$ differentiable in $\ba$, define the policy gradient of $h$ through $\tilde{\bpi}_{\bmu}$ under $\nu$ by
\[
\nabla_{\bmu}^{\mathrm{PG}} \E_{\nu}\!\left[h(\bs, \tilde{\bpi}_{\bmu}(\bs, \blam))\right] \;:=\; \E_{\nu}\!\left[\,\frac{\partial \ba}{\partial \bmu}^{\!\top}\!\nabla_{\ba} h(\bs, \ba)\,\right]_{\ba = \tilde{\bpi}_{\bmu}(\bs, \blam)},
\]
where only the dependence of $\tilde{\bpi}_{\bmu}$ on $\bmu$ is differentiated; any explicit $\bmu$-dependence of $h$ (or of $\nu$) is held fixed. Then the stationarity condition \eqref{eq:moment-condition} is equivalent to
\begin{equation}
\label{eq:pge}
\nabla_{\bmu}^{\mathrm{PG}} \E_{\rho^{\bth} \otimes p}\!\left[Q^{\bth}(\bs, \tilde{\bpi}_{\bmu}(\bs, \blam))\right] \;=\; \nabla_{\bmu}^{\mathrm{PG}} \E_{\rho^{\bth} \otimes p}\!\left[f_{\bmu}(\tilde{\bpi}_{\bmu}(\bs, \blam); \bs)\right].
\end{equation}
\end{theorem}

The equation asserts that, through the dual policy and under $\rho^{\bth}\otimes p$, the gradient in $\bmu$ of the true objective, with the upstream action-value $Q^{\bth}$, equals the gradient obtained by treating the surrogate $f_{\bmu}$ as the cost.

\paragraph{Relation to compatible function approximation.}
Compatible function approximation \citep{sutton1999policy,silver2014deterministic} also requires the critic to be at an optimum of its fitting loss, but there a linear critic, with features given by the Jacobian of the policy being improved, yields that policy's true gradient. Here a single parameter $\bmu$ defines both the critic $f_{\bmu}$ and the dual policy, no compatibility condition is imposed, and the equality holds only for gradients in $\bmu$, in expectation under $\rho^{\bth}\otimes p$. It implies nothing for prices outside the support of $p$; the error of the decisions a downstream problem obtains is controlled instead by Proposition~\ref{prop:decerr}, for prices in the support of $p$.

\section{Training}
\label{sec:algorithm}

\subsection{Procedure}
\label{subsec:procedure}

\paragraph{Rollout mechanics.} Training relies on differentiable rollouts: a length-$T$ trajectory is rolled out against a sampled exogenous trace $\{\bo_t\}$, and the differentiability of the transition $\mathcal{T}$ and cost $c$ in the action allows backpropagation through the whole trajectory. At a perturbation step a dual price $\blam \sim p(\blam)$ is drawn and the action is the perturbed argmin $\ba_t = \nabla f_{\bmu}^*(\blam; \bs_t)$, with per-step loss $c(\bs_t, \ba_t) - \blam^\top \ba_t$; at non-perturbation steps the action follows the frozen base policy $\bpi_{\bth}$ with loss $c(\bs_t, \ba_t)$. The perturbed action produces a counterfactual next-state at which the value target $V^{\bth}$ is evaluated, and $K$ common-random-number draws $\blam_1,\ldots,\blam_K$ per step reduce variance.

\paragraph{Perturbation schedule.} Let $M$ be the number of steps between consecutive perturbations. When the value target is estimated by Monte Carlo rollout continuation, later perturbations at $t+M, t+2M, \ldots$ contaminate the step-$t$ target by $\sum_{j\geq 1} \gamma^{jM}\E_{\blam}[A^{\bth}(\bs_{t+jM}, \ba(\blam,\bs_{t+jM};\bmu))]$. What matters for the gradient is the action sensitivity $\|\partial \bs_{t+jM}/\partial \ba_t\|$, which decays once $jM$ exceeds the closed-loop mixing time $\tau_{\mathrm{mix}}$: a stabilizing policy drives the state back to target after a perturbation, breaking the sensitivity chain (in our lost-sales example an under-order is corrected in one period, while an over-order of size $e$ drains in about $e/\mu$ periods, Appendix~\ref{app:advantage}). The direct dependence of later perturbed actions on $\bmu$ contributes to the gradient at their own steps. The estimator uses $\lfloor T/M\rfloor$ perturbed steps per trajectory with $K$ price draws each; the draws share the state, so under weak dependence along the trajectory the variance scales as $O\big(\tfrac{M}{T}(\sigma_s^2 + \sigma_\lambda^2/K)\big)$, with $\sigma_s^2$ the across-state and $\sigma_\lambda^2$ the within-state variance, while the bias from later perturbations decays once $M$ exceeds $\tau_{\mathrm{mix}}$. The choice $M \approx \tau_{\mathrm{mix}}$ balances the two.

\paragraph{Two-stage procedure.} Stage 1 trains $\bpi_{\bth}$ to minimize the exo-MDP cumulative cost by any method yielding a deterministic policy and an evaluable $V^{\bth}$; we use DirectBackprop/HDPO \citep{madeka2022deep, alvo2023neural}, with $V^{\bth}$ from a critic or rollout continuation. Stage 2 fixes $(\bpi_{\bth}, V^{\bth})$ and optimizes $\bmu$ by SGD on the empirical version of \eqref{eq:training-objective}, backpropagating through the differentiable argmin and transition.

\paragraph{Required inputs.} Besides the frozen policy, Stage 2 needs three inputs. First, the continuation value $V^{\bth}$ at counterfactual next states, from a critic or from rollouts under $\bpi_{\bth}$. In an exo-MDP the transition and cost are known functions of the action, so the next state of any perturbed action is obtained by replaying the stored exogenous trace, with no new interaction with the system. Second, differentiability of $\mathcal{T}$ and $c$ in the action for backpropagation, or else a zeroth-order gradient estimator. Third, the perturbation distribution $p(\blam)$, which by Corollary~\ref{prop:ident} should cover the prices the consumer will impose: historical shadow prices of the coupling constraints, the range of vendor discounts, or, as in Section~\ref{sec:experiments}, the range set by the underage and overage costs.

\paragraph{Policy quality and stochastic policies.} The surrogate targets the advantage of the given policy $\bpi_{\bth}$, whatever its quality. The anchored family additionally assumes that $\bpi_{\bth}(\bs)$ minimizes $Q^{\bth}(\bs,\cdot)$; otherwise the structural condition $\nabla_{\ba} f_{\bmu,\bth}(\ba^*;\bs) = \mathbf 0$ is misspecified by the first-order suboptimality $\nabla_{\ba} Q^{\bth}(\bs,\bpi_{\bth}(\bs))$. For stochastic policies, \eqref{eq:training-objective} needs only evaluations of $Q^{\bth}$, which remain defined, and the perturbed action $\nabla f^*_{\bmu}(\blam;\bs)$ remains deterministic, but the anchored family is misspecified: $\E_{\ba\sim\bpi_{\bth}}A^{\bth}(\bs,\ba) = 0$, so $A^{\bth}$ takes negative values that an anchored $f_{\bmu} \ge 0$ cannot represent, and the best it can match is $Q^{\bth} - \min_{\ba} Q^{\bth}$.

\subsection{Network Architecture and Design Choices}
\label{subsec:architecture}

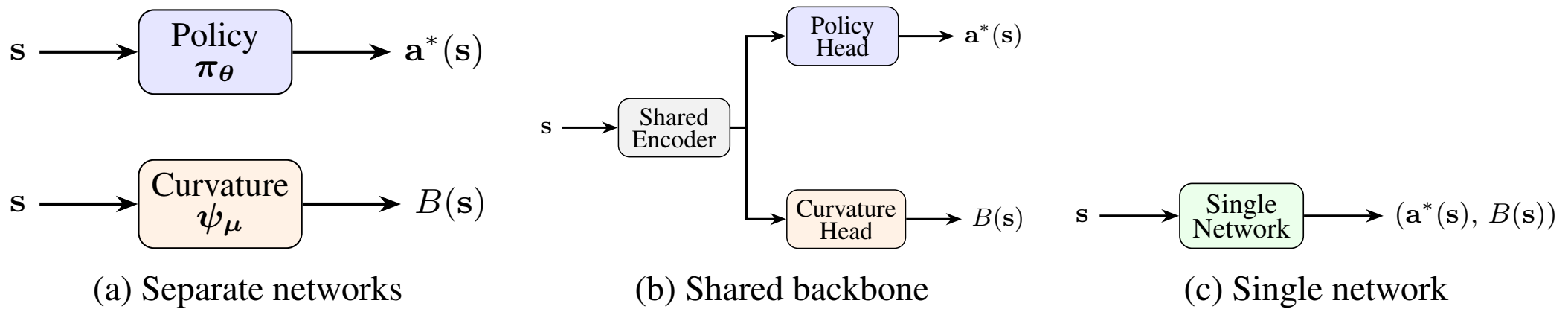
\begin{figure}[htbp]
\centering

\begin{subfigure}[b]{0.32\textwidth}
\centering
\resizebox{\textwidth}{!}{%
\begin{tikzpicture}[
    node distance=0.5cm and 0.9cm,
    block/.style={draw, rounded corners, minimum height=0.7cm, minimum width=1.4cm, align=center, font=\small},
    io/.style={font=\small},
    arr/.style={-{Stealth[length=2mm]}, thick},
]
\node[io] (s1) {$\bs$};
\node[block, right=of s1, fill=blue!10] (pol) {Policy\\[-2pt]$\bpi_{\bth}$};
\node[io, right=of pol] (astar) {$\ba^*(\bs)$};
\draw[arr] (s1) -- (pol);
\draw[arr] (pol) -- (astar);

\node[io, below=1.0cm of s1] (s2) {$\bs$};
\node[block, right=of s2, fill=orange!10] (sur) {Curvature\\[-2pt]$\bpsi_{\bmu}$};
\node[io, right=of sur] (psi) {$B(\bs)$};
\draw[arr] (s2) -- (sur);
\draw[arr] (sur) -- (psi);
\end{tikzpicture}}
\caption{Separate networks}
\end{subfigure}\hfill
\begin{subfigure}[b]{0.32\textwidth}
\centering
\resizebox{\textwidth}{!}{%
\begin{tikzpicture}[
    node distance=0.5cm and 0.7cm,
    block/.style={draw, rounded corners, minimum height=0.7cm, minimum width=1.4cm, align=center, font=\small},
    io/.style={font=\small},
    arr/.style={-{Stealth[length=2mm]}, thick},
]
\node[io] (s) {$\bs$};
\node[block, right=of s, fill=gray!10] (enc) {Shared\\[-2pt]Encoder};
\node[block, above right=0.4cm and 0.7cm of enc, fill=blue!10] (polh) {Policy\\[-2pt]Head};
\node[block, below right=0.4cm and 0.7cm of enc, fill=orange!10] (curh) {Curvature\\[-2pt]Head};
\node[io, right=of polh] (astar) {$\ba^*(\bs)$};
\node[io, right=of curh] (psi) {$B(\bs)$};
\draw[arr] (s) -- (enc);
\draw[arr] (enc.east) -- ++(0.2,0) |- (polh.west);
\draw[arr] (enc.east) -- ++(0.2,0) |- (curh.west);
\draw[arr] (polh) -- (astar);
\draw[arr] (curh) -- (psi);
\end{tikzpicture}}
\caption{Shared backbone}
\end{subfigure}\hfill
\begin{subfigure}[b]{0.32\textwidth}
\centering
\resizebox{\textwidth}{!}{%
\begin{tikzpicture}[
    node distance=0.5cm and 0.9cm,
    block/.style={draw, rounded corners, minimum height=0.7cm, minimum width=1.4cm, align=center, font=\small},
    io/.style={font=\small},
    arr/.style={-{Stealth[length=2mm]}, thick},
]
\node[io] (s) {$\bs$};
\node[block, right=of s, fill=green!8] (net) {Single\\[-2pt]Network};
\node[io, right=of net] (out) {$(\ba^*(\bs),\, B(\bs))$};
\draw[arr] (s) -- (net);
\draw[arr] (net) -- (out);
\end{tikzpicture}}
\caption{Single network}
\end{subfigure}

\caption{Three architecture options for the advantage surrogate. (a)~Separate policy and curvature networks. (b)~Shared encoder with specialized heads. (c)~Single network outputting both decision and curvature parameters.}
\label{fig:architectures}
\end{figure}

The surrogate requires two outputs per state: the anchor $\ba^*(\bs)$ and the curvature parameters $\bpsi(\bs)$ (or $B(\bs)$ for a quadratic). Three architectures are natural (Figure~\ref{fig:architectures}): separate policy and curvature networks, a shared backbone with specialized heads, or a single network outputting both. Separate networks are the most modular option and the natural choice for two-stage training and for retrofitting a surrogate to a pre-trained policy; we use them throughout, with the head architecture of Appendix~\ref{app:exp-arch}.

For the surrogate parameterization, the anchored quadratic $\psi_{\bmu}(\bu; \bs) = \tfrac{1}{2}\bu^\top B(\bs; \bmu)\bu$ is the simplest option, although problems with asymmetric or heavy-tailed advantage geometry admit any strongly convex family (piecewise-quadratic, or a quadratic plus an input convex network \citep{amos2017input}) whose conjugate gradient is evaluable.

\section{Experiments}
\label{sec:experiments} 

\subsection{Setup}
\label{subsec:exp-setup}

We use the exo-MDP formulation of stationary multi-period inventory management with lost sales and zero lead time \citep[Section~4.2]{maggiar2025structure}: endogenous on-hand inventory $y_t$, exogenous demand history $\{d_{t-s}\}_{s=1,\ldots,H}$, static economic parameters (price $p$, cost $c$, holding cost $h$, lost-sale penalty $b$), and order action $a_t \in \R_{\geq 0}$, with $y_{t+1}=(y_t + a_t - d_t)^+$ and $c_t(s,a)=ca + b(d_t-y-a)^+ + h(y+a-d_t)^+-p\min(d_t,y+a)$. The problem admits closed-form optima: the policy is an order-up-to level at a demand quantile,
\begin{align}
a_t^*&= \left(F_t^{-1}\!\left(\frac{c_u}{c_u+c_o}\right) - y\right)^+\label{eq:example_opt_action}\\
\hat G_t(a,s)&= c_o\E\left[(a+y-D_t)^+\right] + c_u \E\left[(D_t - (a+y))^+\right] + K\label{eq:example_q_func}
\end{align}
where $c_u=p-c+b$ and $c_o=h$ are the underage and overage costs (with $\gamma = 1$ and leftover units credited at cost $c$, Appendix~\ref{app:advantage}), $K$ does not depend on $a$, and $F_t^{-1}$ is the inverse cumulative distribution function of the demand distribution $D_t$. The quantity $\hat G_t$ is the single-period newsvendor cost, which is not the advantage of the base-stock policy.

\paragraph{True advantage.} Write $e = a - a^*$ for a one-shot deviation from the optimal order, after which the policy reverts to its base-stock level $z^*$, and $g(e) = \hat G(z^*+e) - \hat G(z^*)$ for the single-period cost of that deviation. An under-order ($e<0$) is corrected by the next order and, with lost sales, leaves no carryover, so its advantage is the single-period cost, $A^-(e) = g(e)$. An over-order ($e>0$) cannot be undone, since orders are nonnegative; the excess is carried until cumulative demand $S_j = \sum_{i<j} D_{t_0+i}$ drains it, and its advantage is the renewal drain
\begin{equation}
\label{eq:advantage-over}
A^{+}(e) = \sum_{j \ge 0} \E\!\big[\, g\big((e - S_j)^{+}\big)\,\big],
\end{equation}
a multi-period cost growing as $\tfrac{c_o}{2\mu}e^2$ for large $e$, where $\mu$ is the mean demand. The true advantage $A^{\bth}$ therefore coincides with $\hat G$ below the optimum and exceeds it above, where the cost of the deviation is carried over several periods. Appendix~\ref{app:advantage} gives the derivation.

We compare two anchored surrogates, motivated by the asymmetry of underage versus overage costs in inventory management: the \textbf{symmetric quadratic} $f_\mu(a; \bs) = \tfrac{1}{2}B(\bs)(a - a^*(\bs))^2$ and the \textbf{asymmetric (piecewise) quadratic} $f_{\bmu}(a; \bs) = \tfrac{1}{2}B_+(\bs)[(a - a^*)^+]^2 + \tfrac{1}{2}B_-(\bs)[(a - a^*)^-]^2$, both with positive curvatures and closed-form conjugates. We report two-stage training at perturbation cadences $M \in \{1, 10, 100\}$.

\paragraph{Training.} We generate $N_{\mathrm{train}}=40{,}000$ products using the data-generating procedure of \citet[Section~4.1 and Table~1]{maggiar2025structure}, with a different initial inventory (32 initialization + 100 rollout periods, Appendix~\ref{app:exp-data}), whose Gamma-distributed demand admit the closed forms \eqref{eq:example_opt_action}--\eqref{eq:example_q_func}. Stage-1 trains the policy network and Stage-2 the curvature network on top of it, both for 500 epochs with separate dilated-CNN policy and curvature heads (full architecture and hyperparameters in Appendix~\ref{app:experiments}). The single-period cost $\hat G$ is asymptotically linear with slopes $-c_u$ and $c_o$, which sets the scale of the dual prices; we draw the training perturbation $\lambda_t \sim \tfrac{1}{2}(c_o-c_u) + \sigma(c_o+c_u)\mathcal{N}(0,1)$ with $\sigma=0.33$ at every $M$-th post-initialization period and $\lambda_t=0$ otherwise, centered and scaled so its bulk falls inside $[-c_u, c_o]$, and cap the surrogate action at $5\,a^*(\bs)$ (Appendix~\ref{app:exp-surr}).

\paragraph{Evaluation.} We evaluate on $N_{\mathrm{test}}=100{,}000$ fresh products. The \emph{unconstrained} evaluation rolls the policy at $\blam \equiv 0$ over 500 periods and regresses the trained order quantities against the non-omniscient closed-form newsvendor optimum \eqref{eq:example_opt_action}, with $F_t$ fit to the demand history the policy sees. The \emph{constrained} evaluation realizes the shared-capacity use case of Section~\ref{sec:intro}: a transient storage-capacity limit on units held, common across products, is dualized to a single scalar shadow price $\lambda$ that surcharges the inventory position held under the constraint. We deliver it as the dual perturbation $\delta = -\lambda$ (a shared negative price), so the surrogate-induced action is the closed-form perturbed argmin $\nabla f_{\hat{\bmu}}^*(-\lambda; \bs)$, which orders less. We evaluate a single constrained period from a null initial inventory (Appendix~\ref{app:capacity}) and report both the realized reward $\overline{R}$ and the mean per-product on-hand position, sweeping $\lambda \in \{0, 1, 2, 4, 8, 16\}$.

\subsection{Results}
\label{subsec:exp-results}

\subsubsection{Training}
\label{subsec:results-training}

\begin{figure}[htbp]
  \centering
  \includegraphics[width=0.8\linewidth,keepaspectratio]{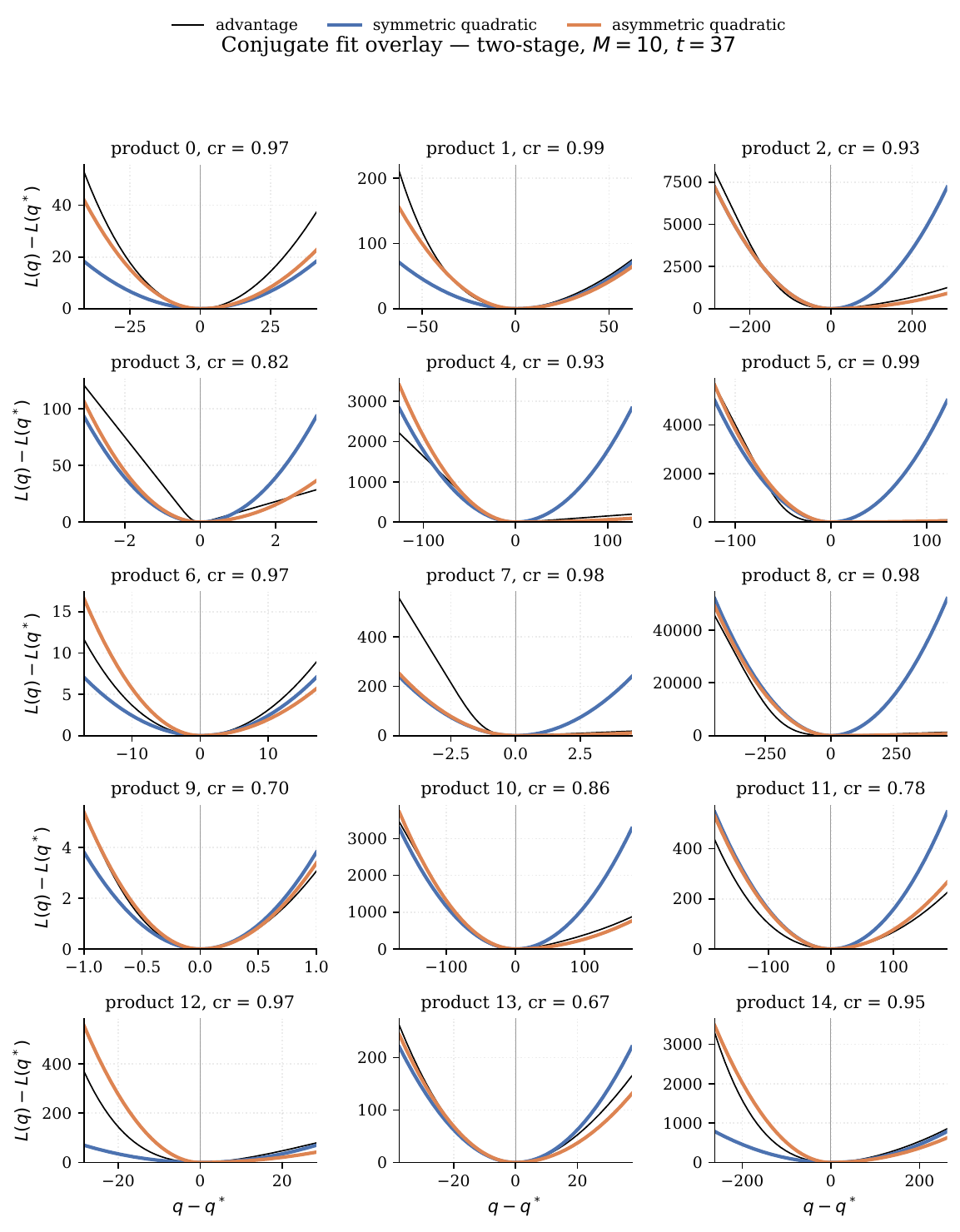}
  \caption{Conjugate fits for the first 15 products of the test set at evaluation step $t = 37$, two-stage training at $M = 10$. Each panel overlays the true one-shot-deviation advantage $A^{\bth}$ (thin black, centered at the anchor; the single-period cost on the underage arm and the multi-period renewal drain on the overage arm, Appendix~\ref{app:advantage}), with the trained symmetric (blue) and asymmetric (orange) quadratic surrogates. Titles report the product's critical ratio.}
  \label{fig:conjugate-fit-overlay}
\end{figure}

The Stage-2 loss is noisier at $M = 1$, where a price is injected at every step and the inventory dynamics never relax, and stabilizes as $M$ grows (Figure~\ref{fig:loss-curves}, Appendix~\ref{app:extra}). Figure~\ref{fig:conjugate-fit-overlay} shows the surrogate fits at $M = 10$: the two surrogates agree near the anchor where the local quadratic is unconstrained, the asymmetric variant is tighter on the underage tail (low critical ratio) on most products, and, as expected of a quadratic, neither follows the linear underage arm of $A^{\bth}$ far from the anchor.

\subsubsection{Unconstrained Evaluation}
\label{subsec:results-uneval}

The unconstrained evaluation reports the mean per-period reward $\overline{R}$ and the slope and $R^2$ of regressing trained order quantities against the closed-form non-omniscient optimum \eqref{eq:example_opt_action} over evaluation steps $t \in [32, 499]$. At $\blam = 0$ the curvature head is inactive, so the action reduces to the frozen Stage-1 anchor $\nabla f^*_{\hat{\bmu}}(0;\bs) = \ba^*(\bs)$ irrespective of the surrogate shape. This evaluation thus checks the anchor alone, before the curvature is tested in the constrained sweep. The anchor reproduces the optimum at slope $0.93$, $R^2 = 0.90$ (Table~\ref{tab:results-unconstrained}), and the trained-vs.-optimal scatter is a tight diagonal band (Figure~\ref{fig:regressions-unconstrained}).

\begin{table}[htbp]
\centering
\begin{minipage}[c]{0.5\linewidth}
\centering
\caption{Unconstrained evaluation ($\blam \equiv 0$). $\overline{R}$ is the mean per-period reward; slope and $R^2$ regress trained order quantities against the closed-form non-omniscient newsvendor optimum over evaluation steps $t \in [32, 499]$, through the origin and excluding zero orders.}
\label{tab:results-unconstrained}
\begin{tabular}{lrrr}
\toprule
Mode & $\overline{R}$ & Slope & $R^2$ \\
\midrule
Stage-1 anchor & $4{,}593.69$ & $0.93$ & $0.90$ \\
\bottomrule
\end{tabular}
\end{minipage}\hfill
\begin{minipage}[c]{0.45\linewidth}
\centering
\includegraphics[width=\linewidth]{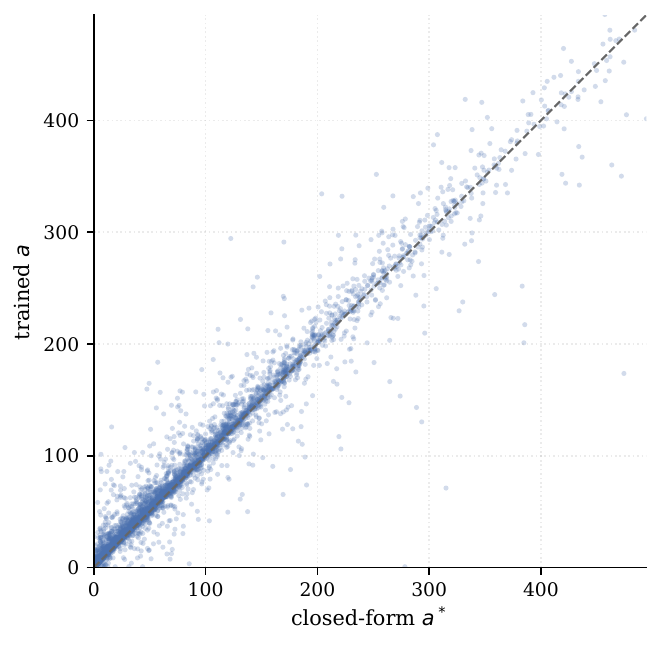}
\captionof{figure}{Trained actions vs.\ the closed-form non-omniscient newsvendor optimum, $\blam \equiv 0$, periods $t \in [100, 109]$.}
\label{fig:regressions-unconstrained}
\end{minipage}
\end{table}

\subsubsection{Constrained Evaluation}
\label{subsec:results-ceval}
The constrained evaluation probes the conjugate response to a shared capacity price: the constrained-period action is $\nabla f_{\hat{\bmu}}^*(-\lambda; \bs) = \ba^*(\bs) + \nabla \psi_{\hat{\bmu}}^*(-\lambda;\bs)$, which trims the order below the anchor by an amount set by the surrogate's curvature. We compare the trained asymmetric quadratic surrogate against two references: a \emph{passive anchor} that ignores the price and plays $\ba^*(\bs)$ regardless (isolating the value of responding to the constraint), and the \emph{non-omniscient optimal}, the exact penalized-newsvendor order-up-to computed from the demand history the policy sees (Appendix~\ref{app:capacity}). Table~\ref{tab:results-capacity} reports the mean per-period realized reward and the mean per-product on-hand position; Figure~\ref{fig:reward-vs-capacity} (Appendix~\ref{app:extra}) plots both against $\lambda$.

As $\lambda$ rises, the passive anchor holds a fixed position and pays a linearly growing surcharge, while the surrogate draws its held inventory down toward the optimal reduction, recovering most of the achievable reward. The unchanged upstream policy thus honors a capacity price it was never trained against.

\begin{table}[htbp]
\centering
\begin{small}
\caption{Shared-capacity evaluation. For each shadow price $\lambda$: mean per-period realized reward $\overline{R}$ and mean per-product on-hand position.}
\label{tab:results-capacity}
\begin{tabular}{lrrrrrr}
\toprule
 & $\lambda = 0$ & $\lambda = 1$ & $\lambda = 2$ & $\lambda = 4$ & $\lambda = 8$ & $\lambda = 16$ \\
\midrule
\multicolumn{7}{l}{\emph{Mean per-period reward} $\overline{R}$}\\
passive anchor        & $4{,}598.3$ & $4{,}438.3$ & $4{,}278.3$ & $3{,}958.4$ & $3{,}318.6$ & $2{,}038.9$ \\
asymmetric quadratic  & $4{,}598.3$ & $4{,}438.4$ & $4{,}286.8$ & $4{,}005.6$ & $3{,}516.9$ & $2{,}752.6$ \\
non-omniscient optimal& $4{,}611.7$ & $4{,}445.5$ & $4{,}293.8$ & $4{,}018.3$ & $3{,}551.7$ & $2{,}842.8$ \\
\midrule
\multicolumn{7}{l}{\emph{Mean on-hand position}}\\
passive anchor        & $160.0$ & $160.0$ & $160.0$ & $160.0$ & $160.0$ & $160.0$ \\
asymmetric quadratic  & $160.0$ & $154.1$ & $148.3$ & $136.9$ & $117.3$ & $89.6$ \\
non-omniscient optimal& $174.7$ & $155.9$ & $144.7$ & $128.4$ & $105.0$ & $75.5$ \\
\bottomrule
\end{tabular}
\end{small}
\end{table}

\section{Conclusion}
\label{sec:conclusion}

Policies trained for per-unit decisions are often consumed by systems that must depart from the optimum, and these systems need the cost of the departure, which policy learning alone does not provide. Fitting a surrogate through its unperturbed argmin cannot identify this cost. Matching convex conjugates under a designed price distribution identifies it on the actions reachable under prices in the support of that distribution, so the choice of distribution is the main modeling decision and should reflect the shadow prices consumers will impose. Because the loss evaluates the action-value function only at the surrogate's own decisions, the surrogate can be fit to a frozen, deployed policy without the consumers' optimal decisions. The main open question is the surrogate family: richer shapes, such as input convex networks or $M^\natural$-convex functions for placement problems, enlarge the set of action-value functions the surrogate can match, but their conjugate gradients are no longer available in closed form, and the cost of differentiating through them at scale is unmeasured.

\bibliographystyle{iclr2027_conference}
\bibliography{main}

\appendix

\section{Experimental Details}
\label{app:experiments}

This appendix documents the data-generating process, network architecture, optimizer, and surrogate parameterizations used in the experiments of Section~\ref{sec:experiments}, in enough detail to reproduce the reported numbers.

\subsection{Data-Generating Process}
\label{app:exp-data}

Each product is a tuple of static economic parameters $(p, c, h, b)$ paired with a Gamma demand process $D_t$. Parameters are sampled independently across products as
\begin{align*}
\bar d &\sim \mathrm{Exp}(\text{mean}=100), &
\mathrm{CV} &\sim \mathrm{U}(0, 1), \\
p &\sim \mathrm{Exp}(\text{mean}=100), &
c \mid p &\sim \mathrm{U}(0, p), \\
h &\sim \mathrm{Exp}(\text{mean}=5), &
b &\sim \mathrm{U}(0, 10).
\end{align*}
Demand at every period is drawn i.i.d.\ from $\mathrm{Gamma}(\alpha, \beta)$ with $\alpha = 1/\mathrm{CV}^2$ and scale $\beta = \bar d \cdot \mathrm{CV}^2$, so that $\E[D_t] = \bar d$ and $\mathrm{Var}(D_t) = \bar d^2 \mathrm{CV}^2$. Prices and costs are rounded to two decimals and demands to one. The static parameters are constant across periods within a product. The training set has $N_{\mathrm{train}} = 40{,}000$ products generated over $32 + 100 = 132$ periods (32 used for initialization, 100 for the rollout). The test set has $N_{\mathrm{test}} = 100{,}000$ products generated over $32 + 500 = 532$ periods.

The endogenous on-hand inventory is initialized at $y_0 \sim 1.5\,\bar{d}_{\text{init}} \cdot \mathrm{U}(0,1)$, where $\bar{d}_{\text{init}}$ is the mean of the first $32$ demand observations of the product, and is then evolved deterministically by $y_{t+1} = (y_t + a_t - d_t)^+$.

\subsection{Network Architecture}
\label{app:exp-arch}

The policy and curvature heads share the same architecture but live in separate networks (no parameter sharing). Each head is a temporal-convolutional encoder followed by a small MLP decoder, applied to per-product features:
\begin{itemize}[leftmargin=*, itemsep=2pt]
\item \emph{Series features.} The history of demand over a backward window $H=32$, fed as a single-channel time series.
\item \emph{Static features.} Five scalars: $(p, c, h, b, y_t)$ broadcast across the convolutional output.
\item \emph{Encoder.} A dilated causal convolutional stack with dilations $\{1, 2, 4, 8, 16\}$ (i.e.\ $2^i$ for $i = 0, \ldots, \log_2 H - 1$), output channel width $8$, two encoding layers of hidden width $32$, embedding dimension $32$.
\item \emph{Decoder.} A linear MLP (no hidden layers) of width $32$, mapping the encoder embedding to the head's output.
\end{itemize}
The policy head outputs a single scalar passed through a softplus to enforce non-negativity, giving the anchor $a^*(\bs;\bth) \geq 0$. The curvature head outputs $d_{\bmu}$ scalars, each passed through $\mathrm{ELU}(\cdot) + 1 + 10^{-4}$, where $d_{\bmu}$ is the surrogate's parameter dimension. The last linear layer of the curvature decoder is initialized with zero weights and bias $0.5$, so that at initialization every curvature output is $\mathrm{ELU}(0.5) + 1 + 10^{-4} \approx 1.5$; the policy decoder uses the default initialization.

\subsection{Surrogate Parameterizations}
\label{app:exp-surr}

Both surrogates are anchored (the surrogate's argmin is the policy network's anchor $a^*(\bs;\bth)$); we describe them as functions of the deviation $x = a - a^*(\bs;\bth)$.
\begin{itemize}[leftmargin=*, itemsep=3pt]
\item \emph{Symmetric quadratic} ($d_{\bmu}=1$). $f_{\bmu}(x) = \tfrac{1}{2} B(\bs)\, x^2$, with conjugate response $\nabla f_{\bmu}^{*}(\lambda) = \lambda / B$.
\item \emph{Asymmetric quadratic} ($d_{\bmu}=2$). $f_{\bmu}(x) = \tfrac{1}{2} B_{+}(\bs)\, [x_+]^2 + \tfrac{1}{2} B_{-}(\bs)\, [x_-]^2$, with curvature head outputs $(B_{+}(\bs), B_{-}(\bs))$ and conjugate response $\nabla f_{\bmu}^{*}(\lambda) = \lambda / B_{+}$ for $\lambda > 0$ and $\lambda / B_{-}$ otherwise.
\end{itemize}
Each curvature scalar is clamped to $\geq 10^{-8}$ before use. The action passed to the cost function is then $a_t = \mathrm{clip}(a^*(\bs;\bth) + \nabla f_{\bmu}^{*}(\lambda_t;\bs);\, 0,\, 5\,a^*(\bs;\bth))$. The upper cap at $5\,a^*$ mirrors the soft truncation discussion of Section~\ref{sec:experiments} and bounds the action's response to dual prices outside $[-c_u, c_o]$.

\subsection{Training Loop}
\label{app:exp-train}

Both stages share the following optimizer and rollout settings:
\begin{itemize}[leftmargin=*, itemsep=2pt]
\item \emph{Optimizer.} AdamW with learning rate $3 \times 10^{-3}$ on the active parameter group, gradient-norm clipping at $1.0$, no weight decay tuning beyond AdamW defaults.
\item \emph{Rollout.} Reset period $32$; rollout length $100$ periods of differentiable backpropagation, with leftover inventory credited at cost through a terminal reward $y_T c$. In Stage 2 the last $10$ rollout periods are unperturbed. Discount factor $\gamma = 1$. Mini-batches of $2{,}000$ products, drawn without replacement.
\item \emph{Epochs.} $500$ epochs for every cell. EMA shadow weights start at epoch $10$ with decay $0.995$; evaluation always uses the EMA snapshot of the final epoch.
\item \emph{Parameter groups.} Stage-1 trains the policy head only ($\lambda \equiv 0$ during the rollout). Stage-2 trains the curvature head only with the Stage-1 policy frozen by exclusion from the optimizer's parameter list.
\end{itemize}
The Stage-1 EMA snapshot is reused as the policy initialization for both Stage-2 surrogates, so the unconstrained policy is held identical across rows in Tables~\ref{tab:results-unconstrained}--\ref{tab:results-capacity}.

\subsection{Shared-Capacity Constrained Example}
\label{app:capacity}

This appendix specifies the shared-capacity constrained evaluation of Section~\ref{subsec:results-ceval} and derives its exact optimal benchmark.

\paragraph{Setting.} A transient storage-capacity limit is common to all products: in the constrained period a bound $\sum_j z_t^j \le C_t$ caps the total inventory held, where $z = y + a$ is the post-order position. Dualizing this coupling constraint attaches a single scalar shadow price $\lambda \ge 0$, identical across products, that surcharges the held position $z$ (the dual term is $-\lambda z$ in the per-period reward). The price is shared, and products respond heterogeneously only through their own economics. We evaluate a single constrained period from a null initial inventory $y = 0$ because the position recovers to the unconstrained base stock as soon as the constraint lifts, thus one period suffices to capture the effect.

\paragraph{Exact optimal order-up-to.} With purchase cost charged at order, inventory carried at holding cost $h$, and terminal salvage at cost, the effective per-unit overage in the single-period newsvendor is $h$ and the underage is $c_u = p - c + b$ (as in \eqref{eq:example_q_func}). Adding the shared surcharge $-\lambda z$ makes the constrained-period objective
\begin{equation*}
\max_{z \ge 0} \big\{\, -c_u\,\E[(D-z)^+] - h\,\E[(z-D)^+] - \lambda z \,\big\},
\end{equation*}
whose first-order condition $c_u - (c_u + h)F(z) - \lambda = 0$ gives the penalized order-up-to
\begin{equation}
\label{eq:capacity-opt}
z_\lambda = F^{-1}\!\left(\left(\frac{c_u - \lambda}{c_u + h}\right)^{+}\right),
\end{equation}
where the critical ratio is clipped at zero: when $\lambda \ge c_u$ the marginal held unit cannot cover the surcharge and the optimum is not to buy ($z_\lambda = 0$). Since $z_\lambda \le z^*$, the constrained period stays on the underage arm, where the single-period cost is the exact advantage (Appendix~\ref{app:advantage}), so the benchmark is exact. The \emph{non-omniscient optimal} of Table~\ref{tab:results-capacity} evaluates \eqref{eq:capacity-opt} with $F$ the Gamma fit to the product's 32-period demand history (the same information the policy sees).

\paragraph{Surrogate response and delivery.} The surrogate's constrained-period action is the anchored perturbed argmin $\ba^*(\bs) + \nabla \psi_{\hat{\bmu}}^*(-\lambda; \bs)$, delivered by injecting the shared dual value $\lambda_t = -\lambda$ at the constrained period (zero otherwise). The negative price drives the response down the underage arm (a change of $-\lambda/B_-$ for the asymmetric quadratic), trimming the order below the anchor. The passive anchor plays $\ba^*(\bs)$ and pays $-\lambda\,z$ on its unadjusted position, isolating the value of responding to the price.

\subsection{True Advantage of the Lost-Sales Base-Stock Policy}
\label{app:advantage}

This appendix derives the true advantage $A^{\bth}(\bs,a) = Q^{\bth}(\bs,a) - V^{\bth}(\bs)$ of the base-stock policy in the lost-sales, zero-lead world of Section~\ref{subsec:exp-setup}, and shows how it differs from the single-period newsvendor cost $\hat G$ of \eqref{eq:example_q_func}. With zero lead time, $\gamma = 1$ and leftover units credited at cost, minimizing $\hat G$ recovers the optimal stationary base-stock level (\citealp{zipkin2000foundations}; \citealp[Section~4.2]{maggiar2025structure}), but the advantage measures something different: the cost of a \emph{one-shot} deviation that perturbs the order-up-to level in a single period and then reverts to the optimum thereafter. As we show below, $A^{\bth}$ therefore coincides with $\hat G$ only where a deviation carries no multi-period consequence. Write the deviation from the optimal order-up-to level as $e = a - a^*$ (equivalently the post-order position is $z^* + e$ with $z^* = F^{-1}(c_u/(c_u+c_o))$), and let $g(e) = \hat G(z^*+e) - \hat G(z^*)$ be the single-period cost deviation, so $g(0)=0$, $g'(0)=0$, with asymptotic slopes $-c_u$ (under) and $c_o$ (over).

\paragraph{Mechanism.} By definition $Q^{\bth}(\bs,a) = c(\bs,a) + \gamma\,\E[V^{\bth}((z-D)^+)]$ with $z = y + a$. The continuation value $V^{\bth}(y')$ is linear in $y'$ for $y' \le z^*$: the policy is already at its optimum, so only the purchase term $-c\,y'$ carried by a lower starting inventory survives, and this linear piece combines with $c(\bs,a)$ to give exactly $\hat G$. For $y' > z^*$ the non-negativity of orders ($a \ge 0$) prevents the policy from ordering down, so it is stuck above target and pays excess holding until demand drains the surplus; $V^{\bth}$ bends strictly convex above the linear extrapolation. A deviation with $z \le z^*$ keeps $y' = (z-D)^+ \le z^*$ in every future period, so it stays on the linear branch and $\hat G$ is exact; a deviation with $z > z^*$ can leave $y' > z^*$, entering the convex branch and adding a multi-period cost on top of $\hat G$.

\paragraph{Underage arm ($e < 0$).} An under-order is corrected by the very next reorder (a single period), and lost sales leave no carryover, so the advantage is exactly the single-period cost,
\begin{equation*}
A^{-}(e) = \hat G(z^*+e) - \hat G(z^*) = g(e),
\end{equation*}
Here $\hat G$ coincides with $A^{\bth}$.

\paragraph{Overage arm ($e > 0$).} The excess $e$ can only be worked off by cumulative demand $S_j = \sum_{i<j} D_{t_0+i}$. Summing the single-period cost deviation $g$ over the excess still outstanding at each future period gives the renewal drain \eqref{eq:advantage-over} of Section~\ref{subsec:exp-setup}, exact under the assumptions above. It is a multi-period cost (about $e/\mu$ periods for $e \gg \mu$, where $\mu$ is the mean demand) that grows super-linearly. By the elementary renewal theorem the renewal function satisfies $U(x)/x \to 1/\mu$ as $x \to \infty$ \citep{feller1971introduction}, which yields $A^{+}(e) \sim \tfrac{c_o}{2\mu}\, e^2$ (a ratio limit). A heuristic two-term approximation, accurate once $e$ exceeds about one period of mean demand, is $A^{+}(e) \approx \beta_{\mathrm{r}}\, g(e) + \tfrac{1}{\mu}\int_0^{e} g(u)\,du$ with $\beta_{\mathrm{r}} = \tfrac{1}{2}(1 + \mathrm{cv}^2)$.

\section{Proofs of Proposition~\ref{prop:decerr} and Corollary~\ref{prop:ident}}
\label{app:ident}

Drop $\bs$ from the notation and write $Q = Q^{\bth}(\bs,\cdot)$, $f = f_{\bmu}(\cdot;\bs)$, $\ba_{\bmu} = \ba(\blam,\bs;\bmu)$ and $\ba_{\blam} = \ba(\blam,\bs)$. We use standard facts on conjugates of closed convex functions \citep[Sections~23--26]{rockafellar1997convex}.

\paragraph{Proposition~\ref{prop:decerr}.} Fix $\blam \in \operatorname{supp} p$ and let $h(\ba) = Q(\ba) - \blam^\top\ba$. Then $h$ is convex with minimizer $\ba_{\blam}$, so $\nabla h(\ba_{\blam}) = \mathbf 0$ and $h(\ba_{\blam}) = -Q^*(\blam)$. The Fenchel-Young gap at the surrogate's decision is therefore
\[
Q(\ba_{\bmu}) + Q^*(\blam) - \blam^\top\ba_{\bmu} = h(\ba_{\bmu}) - h(\ba_{\blam}).
\]
The segment from $\ba_{\blam}$ to $\ba_{\bmu}$ lies in the convex set $\mathcal{C}$, on which $h$ is $m$-strongly convex with $L$-Lipschitz gradient. With $\nabla h(\ba_{\blam}) = \mathbf 0$, this gives $\tfrac{m}{2}\|\ba_{\bmu} - \ba_{\blam}\|^2 \le h(\ba_{\bmu}) - h(\ba_{\blam}) \le \tfrac{L}{2}\|\ba_{\bmu} - \ba_{\blam}\|^2$. Taking expectations under $\rho^{\bth}\otimes p$, and noting that $\Delta(\bmu)$ is the expected gap, proves the claim.

\paragraph{Corollary~\ref{prop:ident}.} By Proposition~\ref{prop:decerr}, $\Delta(\bmu) = 0$ implies $\ba_{\bmu} = \ba_{\blam}$ for $\rho^{\bth}\otimes p$-almost every $(\bs,\blam)$; fix an $\bs$ for which this holds for $p$-almost every $\blam$. Because $f$ is strongly convex with continuous gradient, $f^*$ is differentiable on $\R^d$ with continuous gradient $\nabla f^*(\blam) = \ba_{\bmu}$. Because $Q$ is closed and strictly convex, $Q^*$ is continuously differentiable on $\operatorname{int}\operatorname{dom} Q^*$, with $\nabla Q^*(\blam) = \ba_{\blam}$. The set of $\blam \in \Lambda$ where these continuous maps differ is open and has $p$-measure zero; since $p$ has a positive density on $\Lambda$, it is empty. Integrating along segments of the convex set $\Lambda$ gives $f^* - Q^* \equiv C$ on $\Lambda$. For $\blam \in \Lambda$ and $\ba = \nabla Q^*(\blam) = \nabla f^*(\blam)$, Fenchel's equality gives $Q(\ba) = \blam^\top\ba - Q^*(\blam)$ and $f(\ba) = \blam^\top\ba - f^*(\blam)$, so $f(\ba) - Q(\ba) = -C$ on $\mathcal{A}_\Lambda$.

If $\mathbf 0 \in \Lambda$, the anchor $\ba^*(\bs) = \nabla f^*(\mathbf 0) = \nabla Q^*(\mathbf 0)$ lies in $\mathcal{A}_\Lambda$, and since the anchor is the policy action by assumption, $Q(\ba^*) = Q^{\bth}(\bs,\bpi_{\bth}(\bs)) = V^{\bth}(\bs)$. Anchoring gives $f(\ba^*) = 0$, hence $C = V^{\bth}(\bs)$ and $f = Q - V^{\bth} = A^{\bth}$ on $\mathcal{A}_\Lambda$. \hfill$\square$

\section{Additional Results}
\label{app:extra}

\begin{figure}[htbp]
  \centering
  \includegraphics[width=0.8\linewidth]{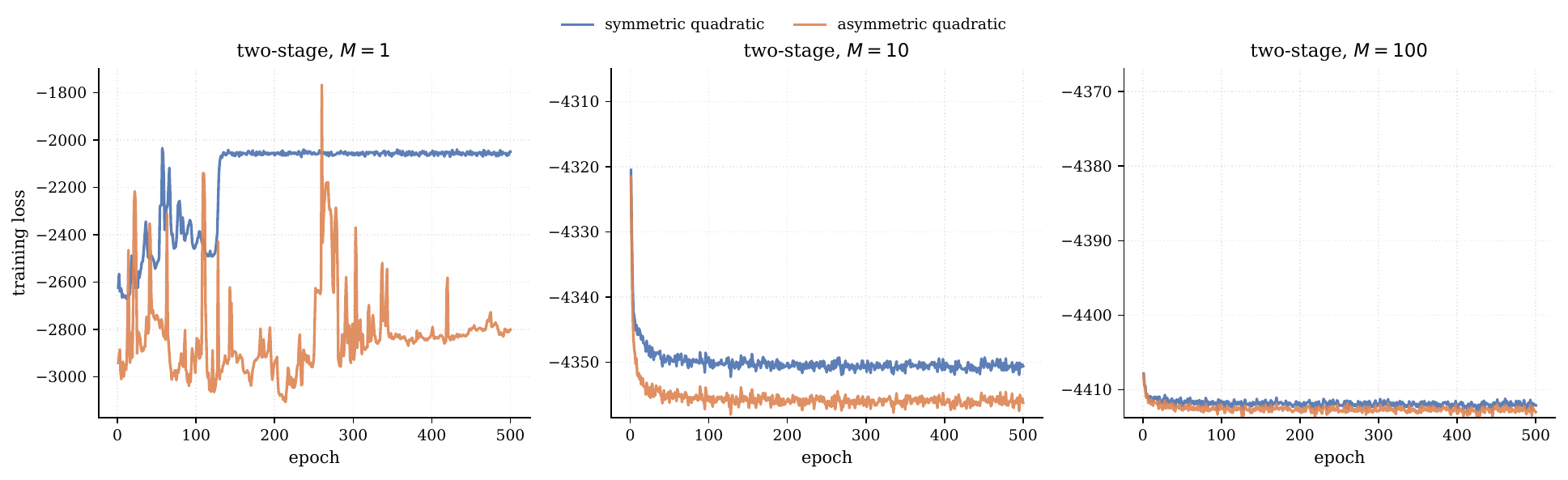}
  \caption{Stage-2 training loss for the symmetric and asymmetric quadratic surrogates, columns $M \in \{1, 10, 100\}$.}
  \label{fig:loss-curves}
\end{figure}

\begin{figure}[htbp]
  \centering
  \includegraphics[width=0.75\linewidth]{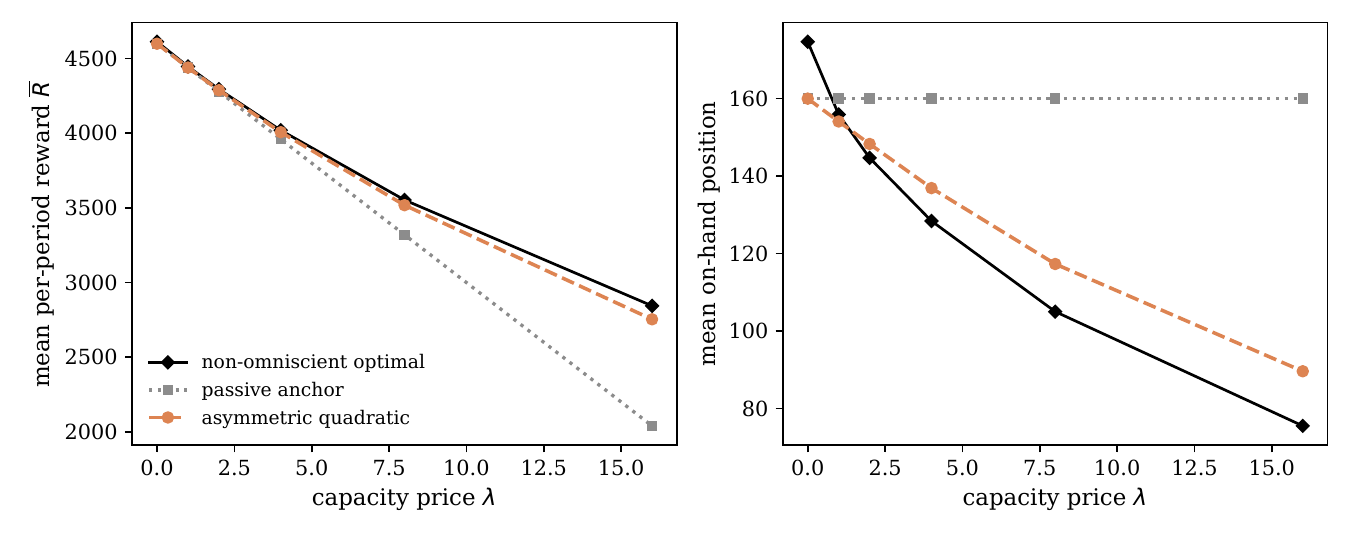}
  \caption{Shared-capacity evaluation vs.\ the capacity price $\lambda$.}
  \label{fig:reward-vs-capacity}
\end{figure}

\end{document}